\documentclass[12pt]{article}%{report}
\usepackage{amsmath,amssymb,array,amsfonts,amsthm}
\usepackage[pdftex,bookmarks,colorlinks=true, anchorcolor=blue, citecolor=black, filecolor=black, menucolor=black, urlcolor=blue, linkcolor=black]{hyperref}
\usepackage{mathrsfs} 
\usepackage[square]{natbib}
\usepackage{bibentry}

\makeatletter 
 \renewcommand{\@oddhead}{\hfill \thepage \hfill}
\makeatother
\renewcommand{\baselinestretch}{1.1} %збільшити інтервал між рядками 
\fussy
\newtheorem{sz_thm}{Theorem}

\newtheorem{sz_col}{Corollary}
\newtheorem{sz_dfn}{Definition}

\newcommand{\nmr}{\symbol{"9D}}

\newcommand{\upR}{\mathrm{R}}

\newcommand{\upF}{\mathrm{F}}
\newcommand{\upH}{\mathrm{H}}
\newcommand{\upC}{\mathrm{C}}

\newcommand{\upE}{\mathrm{E}}

\newcommand{\upQ}{\mathrm{Q}}
\newcommand{\upL}{\mathrm{L}}

\newcommand{\M}{\mathrm{E}}

\newcommand{\upX}{\mathrm{X}}

\newcommand{\bbL}{\mathbb{L}}
\newcommand{\bbW}{\mathbb{W}}

\newcommand{\im}{\mathrm{R}\,}
\newcommand{\tr}{\mathop{\mathrm{tr}}}

\newcommand{\ov}{\mathscr{D}\,}

\newcommand{\ddt}{\mathop{{d\over dt}}}

\newcommand{\dt}{\mathop{\mathrm{dt}}}

\renewcommand{\le}{\leqslant}

\newcommand{\tilt}{\Leftrightarrow}

\newcommand{\dfn}{\mathop{\stackrel{\mathrm{def}}{=}}}

\newcommand{\inpo}[1][\cdot,\cdot]{%
     \def\arg1{{#1}}%
     \spRelay%
}
\newcommand{\spRelay}[1][]{%
  \ensuremath{
    \mathop{\langle}\arg1\rangle_{\ssstyle #1}%
  }
}
\newcommand{\ip}[2]{\inpo[#1][#2]}

\newcommand{\ipll}[1]{\inpo[#1][\upL_2^l]}

\newcommand{\iprm}[1]{(#1)_{\ssstyle m}}
\newcommand{\iprn}[1]{(#1)_{\ssstyle n}}
\newcommand{\iprl}[1]{(#1)_{\ssstyle l}}
\newcommand{\iprp}[1]{(#1)_{\ssstyle p}}

\newcommand{\nrm}[1][\cdot]{%
     \def\arg1{{#1}}%
     \nrmRelay%
}
\newcommand{\nrmRelay}[1][]{%
  \ensuremath{
    \bigl\|\arg1\bigr\|^2_{\ssstyle #1}%
  }
}

\newcommand{\nrLm}[1]{\nrm[#1][\upL^m_2]}
\newcommand{\nrLl}[1]{\nrm[#1][\upL^l_2]}

\newcommand{\Fx}{\upF\x}
\newcommand{\Fz}{\upF'\z}

\newcommand{\oH}{\ensuremath{\mathcal{H}}}

\newcommand{\oHs}{{\mathcal{H}^*}}
\newcommand{\oT}{\ensuremath{\mathcal{T}}}

\newcommand{\oD}{\ensuremath{\mathcal{D}}}

\newcommand{\oDs}{\ensuremath{\mathcal{D}^*}}

\newcommand{\imT}{\im(\oT)}
\newcommand{\setU}{\ensuremath{\mathscr{U}}}

\newcommand{\setW}{\ensuremath{\mathscr{W}}}

\newcommand{\setL}{\ensuremath{\mathscr{F}}}

\newcommand{\Rn}{\ensuremath{\mathbb{R}^n}}

\newcommand{\Rm}{\mathbb{R}^m}
\newcommand{\Rl}{\mathbb{R}^l}
\newcommand{\Rp}{\mathbb{R}^p}

\newcommand{\Wtwone}[3]{
        \ensuremath{
                \bbW^{\ssstyle#1}_{\ssstyle #2}\bigl([a,c],#3\bigr)
        }
}

\newcommand{\Wm}{\Wtwone{1}{2}{\Rm}}

\newcommand{\WnF}{\Wtwone{1}{2,\upF}{\Rn}}

\newcommand{\ssstyle}{\scriptscriptstyle}

\newcommand{\Ltwo}[1]{
        \ensuremath{
                \bbL_{\ssstyle 2}\bigl([a,c],#1\bigr)
        }
}

\newcommand{\Ln}{\Ltwo{\Rn}}

\newcommand{\Lm}{\Ltwo{\Rm}}
\newcommand{\Ll}{\Ltwo{\Rl}}

\newcommand{\eR}{\mathrm{R}_\eta}

\newcommand{\Qw}{\mathrm{Q}_{\ssstyle 2}}

\newcommand{\iQw}{\mathrm{Q}_{\ssstyle 2}^{\ssstyle-1}}

\newcommand{\iQn}{\ensuremath{\mathrm{Q}_{\ssstyle 1}^{\ssstyle-1}}}

\newcommand{\Qn}{\mathrm{Q}_{\ssstyle 1}}
\newcommand{\Qz}{\mathrm{Q}_0}

\newcommand{\iQz}{\mathrm{Q}_0^{\ssstyle -1}}

 \newcommand{\Ktt}{\mathrm{K}'(t)}

\newcommand{\mA}{\mathrm{F}}
\newcommand{\mB}{\mathrm{B}}
\newcommand{\mH}{\mathrm{H}}

\newcommand{\Ct}{\mathrm{C}(t)}

\newcommand{\Kt}{\mathrm{K}(t)}
\newcommand{\Ctt}{\mathrm{C}'(t)}
\newcommand{\Htt}{\mathrm{H}'(t)}
\newcommand{\Ht}{\mathrm{H}(t)}

\newcommand{\tf}{\tilde{f}}

\newcommand{\et}{\eta(t)}
\newcommand{\x}{\ensuremath{x}}
\renewcommand{\v}{\mathrm{v}}
\newcommand{\y}{\ensuremath{\mathrm{y}}}

\newcommand{\hw}{\ensuremath{\hat{\omega}}}

\renewcommand{\u}{\ensuremath{\mathrm{u}}}

\newcommand{\xt}{\x\,(t)}

\newcommand{\pt}{p(t)}

\newcommand{\ut}{u(t)}
\newcommand{\ft}{f(t)}
\newcommand{\yt}{y(t)}

\newcommand{\z}{\mathrm{z}}
\newcommand{\p}{\mathrm{p}}
\newcommand{\zt}{\z\,(t)}

\newcommand{\hp}{\hat{p}}

\newcommand{\hx}{\hat{x}}

\newcommand{\hz}{\hat{z}}

\newcommand{\hu}{\hat{u}}

\newcommand{\hxt}{\hat{x}(t)}
\newcommand{\hzt}{\hat{z}(t)}

 \newcommand{\fz}{\ensuremath{\f_0}}
\newcommand{\lx}{\ensuremath{\ell(\x)}}
\newcommand{\hlx}{\widehat{\ell(\x)}}
\newcommand{\uy}{\ensuremath{\mathrm{u}_c(\y)}}
\newcommand{\huy}{\hu_{\hat{c}}(\y)}

\newcommand{\G}{\mathscr{G}}

\newcommand{\Gw}{\G_{\ssstyle 2}}

\newcommand{\f}{\ensuremath{f}}

\newcommand{\gz}{x^g_0}
\newcommand{\yk}{\y_k}
\newcommand{\pk}{p_k}
\newcommand{\hpk}{\hat{p}_k}
\newcommand{\pko}{p_{k+1}}

\newcommand{\lk}{\ell_k}
\newcommand{\hzk}{\hz_k}

\newcommand{\hzko}{\hz_{k+1}}
\newcommand{\Hk}{\upH_k}
\newcommand{\Hkt}{\upH'_k}

\newcommand{\xk}{\x_k}
\newcommand{\hxk}{\hx_k}

\newcommand{\vk}{\v_k}
\newcommand{\Fk}{\upF_k}
\newcommand{\Fkt}{\upF'_k}
\newcommand{\Fko}{\upF_{k+1}}
\newcommand{\Ck}{\upC_k}
\newcommand{\Ckt}{\upC'_k}

\newcommand{\R}{\upR}

\newcommand{\Qko}{\upQ_{\ssstyle 1,k}}

\newcommand{\Qkw}{\upQ_{\ssstyle 2,k}}

\newcommand{\mybeglist}{\begin{description}\setlength{\itemsep}{.2cm}
}
\newcommand{\myendlist}{
\end{description}}
\newcommand{\mybibentry}[1]{\item[\rm\citep{#1}] \bibentry{#1}}

\begin{document}
\bibliographystyle{plainnat}
%\nobibliography{refs/myrefs,refs/refs}
%\selectlanguage{ukrainian}
%Параметри сторінки
\fontsize{14}{\baselineskip}%встановити шрифт (14pt) і міжрядковий інтервал 
\selectfont %завантажити шрифт
\setlength{\textheight}{36\baselineskip}%встановити кількість рядків у 
%документі рівною Round(36\baselinestretch) у даному випадку 40
\setlength{\textheight}{\baselinestretch\textheight}
\addtolength{\textheight}{\topskip}

\voffset= -4.5mm %зсув цілої сторінки по вертикалі
\hoffset= 5mm %зсув цілої сторінки по горизонталі
\begin{titlepage}\pagestyle{empty}
  \begin{center}
    {Taras Shevchenko National University Of Kyiv}
  \end{center}\vspace{1.5cm}
  \begin{center}
    {\bf\Large Sergiy Zhuk}
  \end{center}\medskip
  \begin{flushright}
    \textbf{UDC 517.926:681.518.2}
  \end{flushright}\vspace{2.5cm}%  \bigskip\bigskip
  \begin{center}
    {\bf Minimax state estimation for linear descriptor systems} 
  \end{center}\vspace{1.5cm}
  \begin{center}%\bigskip
    01.05.04 -- system analysis and optimal decision theory
  \end{center}\vspace{2.5cm}
  \begin{center}
    {Author's Summary}\\
    of the dissertation for the degree of the\\
    Candidate of Science (physics and mathematics)\\
  \end{center}\vspace{2.5cm}
  \begin{center}
    Kyiv-2006
  \end{center}
\newpage
\noindent Dissertation is a manuscript.\\[5mm]
The dissertation accomplished at the 
Department of System Analysis and Decision Making Theory at National Taras Shevchenko University of Kyiv.\\[10mm]
Scientific advisor\\
Doctor of Science (physics and mathematics), Professor Alexander NAKONECHNY, 
Head of the Department of System Analysis and Decision Making Theory at National Taras Shevchenko University of Kyiv.\\[10mm]{
 \emergencystretch=1.5em

 }
\noindent Jury:\\
Doctor of Science (physics and mathematics), Professor Arkadii CHIKRII, 
Glushkov Institute for Cybernetics at National Academy of Science of Ukraine, 
Associate member of the National Academy of Science of Ukraine, Head of the Department of Dynamics of Controlled Processes\\[5mm] 
Doctor of Science (physics and mathematics), Professor Valentin OSTAPENKO, 
Institute for system analysis at Kyiv Polytechnic University, 
Head of the Department of Numerical Methods for Optimization;\\[5mm]
Candidate of science (physics and mathematics), associate professor Alexander 
SLABOSPITSKY,  Department of Applied Statistics at National Taras Shevchenko 
University of Kyiv.\\[10mm]{
 \emergencystretch=1.5em

 }

The manuscript is available at the scientific library of National Taras Shevchenko 
University of Kyiv, Kyiv, 58, Volodymyrska Str.\\[5mm]
\end{titlepage}
{\bf Outline.} This document is organized as follows. Section 1 presents an extended
abstract of the dissertation containing research context, description of the aims of the
dissertation, discussion of the new results, their validation, impact and dissemination,
and information about financial support. Section 2 contains notation, brief summary
of the chapters of the dissertation (definitions, duality concepts and minimax state
estimation algorithms), conclusions, list of papers and conference presentations related
to the dissertation and list of papers cited in the text.

\section{Extended abstract of the dissertation}
  \label{sec:1}
\textbf{Research context.} Differential-Algebraic equations (DAE) are widely used in engineering. In vehicle dynamics, DAEs represent a convenient and powerful basis for software modeling platforms (see for instance the platform \href{www.modelica.org}{Modelica}). DAEs are also used in mathematical economics, robotics, biomechanics, image processing and control theory. In this context, the analysis of DAE is an important issue in the framework of mathematical systems and control theory.  
% as well as in differential equations and numerical analysis. 
%Significant results in these fields of DAE theory were obtained by Gantmaher F., Bojarintsev Y., Kurina G., Rutkas A., Shlapak Y., 
%Luis F., Campbell S., Brenan C., Petzold L., Dai L. and Luenberger D.     
One of key problems of mathematical systems theory is so-called state estimation, that is to construct an estimate of the state, given observations of the state of the process being modeled. The dissertation is devoted to the development of the mathematical theory and algorithms for state estimation problems for linear DAE. 
The main mathematical tool is the Minimax State Estimation (MSE) approach.  
Main assumptions behind MSE are 1) the model is represented by a system of differential equations (for instance, system of Ordinary Differential Equations (ODE) or Partial Differential Equations (PDE)), 2) links between the model and observed data are represented by observation equation and 3) uncertain parameters (for instance, error in initial condition or model error or noise in the observed data) belong to a given \emph{bounding set}. In other words, the uncertainties in the model and observed data are described in terms of the bounding set. The main idea behind MSE is to describe \emph{how the model propagates uncertain parameters which are consistent with observed data and belong to the given bounding set.} 

Main notions of MSE are reachability set, minimax estimate and worst-case error. By definition, reachability set contains all states of the model which are consistent with observed data and uncertainty description. Given a point $P$ within the reachability set one defines a worst-case error as the maximal distance between $P$ and other points of the reachability set. Then the minimax estimate of the state is defined as a point minimizing the worst-case error (a Tchebysheff center of the RS). Basics of the MSE were developed by 
~\cite{Bertsekas1971}, \cite{Tempo1985}, \cite{Kuntsevich1992}, 
\cite{Chernousko1994}, \cite{Kurzhanski1997}, 
\cite{Nakonechny2004}.

In the case of linear model, given a bounding set the classical MSE allows to construct the minimax estimate for the state of the model and calculate the worst-case estimation error, provided the model operator is bounded and has a bounded inverse. 
%The worst-case estimation error $\sigma(u):=\sup_{f,\eta}d(u(y),\varphi)$ 
%depends only on the model $L$, observation operator $H$ and bounding set $
%\mathscr G$. 
The classical MSE is based on the Kalman Duality principle which states that the state estimation problem is equivalent to a dual optimization problem, provided the model operator is bounded and has a bounded inverse. However, linear DAEs do not fit 
this framework as the corresponding model operator may not be invertible or may have unbounded inverse. Therefore, the classical duality concept was not applicable for 
derivation of the MSE theory for linear DAEs. 

\textbf{Aim of the dissertation.} 
The aim of the dissertation is to develop a generalized Kalman Duality 
concept applicable for linear unbounded non-invertible operators and 
 introduce the minimax state estimation theory and algorithms for linear differential-algebraic equations. 
In particular, the dissertation pursues the following goals:
\begin{itemize}
\item develop generalized duality concept for the minimax state estimation theory for DAEs with unknown but bounded model error and random observation noise with unknown but bounded correlation operator;
\item derive the minimax state estimation theory for linear DAEs with unknown but
bounded model error and random observation noise with unknown but bounded
correlation operator;
\item describe how the DAE model propagates uncertain parameters;
\item estimate the worst-case error;  
\item construct fast estimation algorithms in the form of filters;  
\item develop a tool for model validation, that is to assess how good the model describes observed phenomena.
\end{itemize}

\textbf{New results.} The dissertation contains the following new results:
\begin{itemize}
\item generalized version of the Kalman duality principle is proposed 
allowing to handle unbounded linear model operators with 
non-trivial null-space; 
\item new definitions of the minimax estimates for DAEs based on the generalized Kalman duality principle are proposed;
\item theorems of existence for minimax estimates are proved;
\item new minimax state estimation algorithms (in the form of filter and in the variational form) for DAE are proposed. 
\end{itemize}

\textbf{Validation and impact. } 
In order to validate the Generalized Kalman Duality (GKD) concept I applied it~\cite{Zhuk2005a,Zhuk2006b} to linear incorrect differential operators (see \cite{Tikhonov1977} for further details on incorrect problems). As a consequence, I constructed new minimax state estimation algorithms for linear DAE and linear Boundary Value Problems (BVP) for ODE. The main impact is that the new minimax estimate for DAE does not require regularity assumptions on DAE structure (regularity of the matrix pencil~\cite{Hanke1989} or rank-degree condition~\cite{Dai1989}), imposed by the majority of authors; the minimax estimate for BVP works without restricting the structure of the corresponding ODE matrices. In~\cite{Zhuk2005,Zhuk2006b} I applied GKD to the case of the linear discrete-time DAE. As a final result, I obtained a new minimax recursive estimator for discrete-time DAE~\cite{Zhuk2004a,Zhuk2005,Zhuk2005b}. These results demonstrate the impact of the GKD for DAE models. In general, GKD extends the scope of the classical MSE framework, 
originally designed for linear models with bounded invertible operator 
and bounded uncertainties, on the linear models with 
unbounded non-invertible operator and unbounded uncertainties. This brings the following advantages: 1) possibility to construct minimax state estimation algorithms for 
Differential-Algebraic Equations (DAE) and 2) possibility to address the case of unbounded model errors.

\textbf{Dissemination of the new results.}
The generalized Kalman duality concept and minimax state estimation algorithms for linear DAEs were published in the sequence of papers~\cite{Zhuk2004a,Zhuk2005,Zhuk2005b, Zhuk2005a,Zhuk2005b} and was reported at the conferences~\cite{Zhuk2004c,Zhuk2006,Zhuk2006a}. Also the results were presented at the following seminars:  
``System analysis and decision making theory'' and ``Modeling and optimization of uncertain systems'' at the National Taras Shevchenko University of Kyiv, ``Non-smooth optimization'' at the Institute for system analysis at Kyiv Polytechnic University, ``Optimization of controlled processes'' at Glushkov Institute for Cybernetics at National Academy of Science of Ukraine. 

\textbf{Connections with national research programmes.} 
The dissertation was funded by the National State Research Programme \nmr01SF015-01 "Development of the theory, algorithms and software for stochastic and algebraic systems  with applications in economics, engineering and education"  (state registration number is 0101U002173).

\section{Summary of the dissertation}
\label{sec:2}
\textbf{Notation.} 
$E$ denotes the expected value of the random variable,\\
$\mathbb R^n$ denotes the $n$-dimensional Euclidean space,\\
$(\cdot,\cdot)_n$ denotes the canonical inner product in $\mathbb R^n$,\\  
$\mathrm{tr} Q$ denotes the trace of the matrix $Q$,\\
$y=(y_0\dots y_N)$ denotes a vector composed by elements $y_1$,\dots, $y_N$, where 
$y_i$ also may be a vector,\\ 
$\mathbb L_2([a,c],\mathbb R^n)$ denotes a space of 
square-integrable vector-functions on $[a,c]$ with values in $\mathbb R^n$, \\
$\mathbb W_2^1([a,c],\mathbb R^n)$ denotes a space of absolutely continuous vector functions on $[a,c]$ with values in $\mathbb R^n$,\\
$\mathcal H^*$ denotes the adjoint of the linear mapping $\mathcal H$,\\
$B'$ denotes the transposed of the matrix $B$,\\ 
$\mathscr D(L)$ denotes the domain of the linear operator $L$,\\
$\|\cdot\|_{X}$ denotes the norm of the normed space $X$,\\
$\langle\cdot,\cdot\rangle$ denotes the csnonical inner product in the Hilbert space $H$.   

{\bf Outline of the dissertation}. The dissertation is composed of three chapters. The total number of pages is 140. 
The first chapter contains the general description of the dissertation, its aims and basic notions. Also it contains the description of the state of the art in the DAE literature, a brief overview of the state estimation methods and description of the new results. The second chapter describes the generalized Kalman duality concept and minimax state estimation algorithms for linear DAE with discrete time. The third chapter contains the generalized Kalman duality concept and minimax state estimation algorithms for linear DAE with continuous time and conclusions.

{\bf Brief summary.} Let us consider the contents of the second chapter in brief. 
Let $x\in\Rn$ solve a linear algebraic equation 
\begin{equation}
  \label{eq:Fx:a}
  Fx=Bf,
\end{equation} and let the observed data $y\in\Rl$ verify  
\begin{equation} \label{eq:y=Hx+e}
  y=\upH x+\eta,
\end{equation}
where $F$ is $m\times n$-matrix, $\mB$ is $m\times p$-matrix,
 $\mH$ is $l\times n$-matrix and $\eta$ is a realization of the random $l$-vector, 
$x$ denotes the state of the system~\eqref{eq:Fx:a},  
$f\in\Rp$ represents an uncertain element. 

We will assume that $f$ and $
\eR\dfn\upE\eta\eta'$ are uncertain and $$
f\in\G,\eR\in\Gw,
$$ where $\G,\Gw$ are given subsets. In what follows we will be looking for the estimate of the linear function 
$x\mapsto\lx=\iprn{\ell,x},\ell\in\setL$. We will look for the estimate in the class of affine functions  $
y\mapsto\iprl{u,y}+c$ of observed data~\eqref{eq:y=Hx+e}. We will refer $\iprl{u,y}+c$ as an estimate. Let us assign to each estimate $u_c$ an estimation error $$
\sigma(u,c)\dfn\sup_{x,\eR}\{\upE[\iprn{\ell,x}-\iprl{u,y}-c]^2|Fx\in B(\G),
    \eR\in\Gw\},
    $$ where $B(\G)=\{Bf|f\in\G\}$. 
\begin{sz_dfn}
  The estimate $\widehat{(\ell,x)}=\iprl{\hu,y}+\hat{c}$ verifying 
    \begin{equation*}%\label{df:err:a}
      \sigma(\hu,\hat{c})=
      \inf_{u,c}\sup_{x,\eR}\{\upE[\iprn{\ell,x}-\iprl{u,y}-c]^2|Fx\in B(\G),
    \eR\in\Gw\}%=\inf_{u,c}\sigma(u,c)
  \end{equation*}
is called a minimax a priori mean-squared estimate (a priori estimate). 
The number $   \hat{\sigma}=\inf_{u,c}\sigma(u,c)$
   is called a minimax mean-squared error (a priori error).
  \end{sz_dfn}
Let us consider an a posteriori state estimation method. Let $x\in\Rn$ verify~\eqref{eq:Fx:a} and $y\in\Rl$ is given in the form 
\begin{equation}
  \label{eq:y=Hx+g}
  y=Hx+g,
\end{equation}
where $g\in\Rm$ is a vector. In contrast to the previous considerations, we assume that $(f,g)$ are deterministic and belong to the given set $\G\subset\Rl\times\Rm$. Define $$
\upX\dfn\{x\in\Rn|\exists\,(f,g)\in\G:Fx=Bf,y-Hx=g\}
$$ 
\begin{sz_dfn}
  The estimate $\widehat{\ell(x)}$ verifying $$
  \sup_{x\in\upX}|\ell(x)-\widehat{\ell(x)}|=
  \inf_{\tilde{x}\in\upX}\sup_{x\in\upX}|\ell(x)-\ell(\tilde{x})|
  $$ is called a minimax a posteriori estimate (a posteriori estimate). 
  The number $$
  \hat{\sigma}\dfn\sup_{x\in\upX}|\ell(x)-\widehat{\ell(x)}|
  $$ is called minimax a posteriori estimation error (a posteriori error).
\end{sz_dfn}
The a priori and a posteriori estimates are constructed in the dissertation for the
generic convex compact bounding sets $\mathscr G$, $\mathscr G_2$. 
Let us consider the case of ellipsoidal bounding set in more details. 
% Зупинимось більш докладно на випадку квадратичних обмежень.
\begin{sz_thm}\label{tm:aprerrq:a}
  Let $$
  \G=\{f\in\Rp|\iprp{\Qn f,f}\le 1\},\quad\Gw=\{\eR:\tr\Qw\eR\le 1\},
  $$ where $\Qn,\Qw$ are positive definite symmetric matrices. Then the minimax a priori estimate of the linear function $$
  x\mapsto\iprn{\ell,x}, \ell\in\setL\dfn
  \{\ell=F'z+H'u,z\in\Rm,u\in\Rl\}$$ of the solution of $Fx=Bf$ has the following form $$
  \widehat{(\ell,x)}=\iprl{\hu,y},\,\hu=\Qw Hp,\,  
  $$ where $p$ solves 
  \begin{equation}\label{eq:lhzp:a}
    \begin{split}
      &Fp=B\iQn B'\hz,\\
      &F'\hz=\ell-H'\Qw Hp.
    \end{split}
  \end{equation}
The minimax a priori error is given by $$
\sup_{x,\eR}\upE[\iprn{\ell,x}-\widehat{(\ell,x)}]^2=\iprn{\ell,p}
$$
If $\ell\notin\setL$, then the minimax a priori error is infinite. 
\end{sz_thm}
\begin{sz_thm}\label{tm:erraposq:a}
  If $\ell\in\setL$ and 
  \begin{equation*}%  \label{eq:G}
  \G=\{(f,g):\iprp{\Qn f,f}+\iprl{\Qw g,g}\le1\},
\end{equation*}
%%де $\Qn,\,\Qw$ -- додатно означені симетричні матриці, 
then the minimax a posteriori estimate is given by $$
  \hlx=\iprn{\ell,\hx}=\iprl{\Qw Hp,y},
  \hat{\sigma}= [1-\iprl{y-H\hx,\Qw y}]^\frac 12\iprn{\ell,p}^\frac 12
$$ %% Похибка оцінювання дається виразом 
and the minimax a posteriori error is given by $$
     \hat{\sigma}= [1-\iprl{y-H\hx,\Qw y}]^\frac 12\iprn{\ell,p}^\frac 12,$$ 
  where $p$ solves~\eqref{eq:lhzp:a} and $\hx$ solves
  \begin{equation*}%\label{eq:hxhpy}
   \begin{split}
    &F\hx=B\iQn B'\hp,\\
    &F'\hp=H'\Qw(y-H\hx).
  \end{split} 
  \end{equation*}
If $\ell\notin\setL$, then the minimax a posteriori error is infinite.
\end{sz_thm}

Let us demonstrate one application of Theorem~\ref{tm:erraposq:a} to the state estimation for the linear DAEs with discrete time. Assume $x_0\dots x_N$ is a solution of the DAE with discrete time: 
\begin{equation}  \label{eq:Fxk}
  F_{k+1}x_{k+1}-C_k x_k=B_k f_k,\quad F_0\x_0=S\gz,k=\overline{0,N}
\end{equation}
and the observed data is given by
\begin{equation}\label{eq:yk:dae}
  \yk=H_k x_k+g_k,\,k=\overline{0,N}
\end{equation}
where $F_k,C_k,S$ are $m\times n$-matrices, $B_k$ is a $m\times p$-matrix, $f_k\in\Rp$, $\gz\in\Rm$ are some vectors, $
{H_k}$ is $l\times n$-matrix and $g_k\in\Rl$ stands for a deterministic noise in the observed data. Define a linear function
\begin{equation*}%  \label{df:lx:r}
\lx\dfn\sum_{k=0}^{N+1}\iprn{\lk, x_{\ssstyle k}}\dfn(\ell,\x),
\,\lk\in\Rn
\end{equation*} where $\ell:=(\ell_1,\dots,\ell_{N+1})$ and $\x:=(x_0\dots x_N)$. Define $$
F=\left[
\begin{smallmatrix}
  F_0&&0&&0&&\dots&&0&&0\\
  0&&-C_0&&F_1&&\dots&&0&&0\\
  \dots&&\dots&&\dots&&\dots&&\dots&&\dots\\
  0&&0&&0&&\dots&&-C_{N-1}&&F_N
\end{smallmatrix}
\right], H=\mathrm{diag}\{H_0\dots H_N\}, B=\mathrm{diag}\{B_0\dots B_{N-1}\}
$$ and set $y:=(y_0\dots y_N)$, $f:=(\gz,f_0\dots f_{N-1})$, $g:=(g_0\dots g_N)$. It is easy to see that $Fx=Bf$ is equivalent to~\eqref{eq:Fxk} and $y=Hx+g$ is equivalent to~\eqref{eq:yk:dae}. Now we apply Theorem~\ref{tm:erraposq:a} in order to derive the minimax a posteriori estimate of the linear function $(\ell,\x)$ in the variational form.   
\begin{sz_thm}
Assume that $\ell=(\ell_1,\dots,\ell_{N+1})$ is such that $F'_{N+1}z_{N+1}=\ell_{ N+1}$ and $$
F'_k z_k-C'_{k} z_{k+1}+H_{k} u_{k}=\ell_{k},k=\overline{0,N}\,
$$ for some $z_k$ and $u_k$, $k=\overline{0,N}$. 
Assume also that $$
\G\dfn\{(\f,\v):\iprm{\Qz\gz,\gz}+\sum_{k=0}^N\iprm{\Qko f_k,f_k}+
\iprl{\Qkw \vk,\vk}\le 1\},
$$ where $\Qz,\Qko,\Qkw$ are positive definite symmetric matrices for $k\in[0,N]$. Let $(\hpk,\hxk)$ solve  
\begin{equation*}%\label{eq:hxkhpk:r}
  \begin{split}
    &\Fko\hx_{k+1}=\Ck\hxk+B_k\upQ_{\ssstyle 1,k}^{\ssstyle -1}B'_k\hpk,
    \quad F_0\hx_0=S\iQz S'\hp_0,\\
    &\Fkt\hpk=\Ckt\hp_{k+1}+\Hkt\Qkw(\yk-\Hk\hxk),\quad
    {F}'_{\ssstyle N+1}\hp_{\ssstyle N+1}=0,\,k=\overline{0,N}
  \end{split}
\end{equation*}
The minimax a posteriori estimate of the linear function $
\sum_{k=0}^{N+1}\iprn{\lk,\xk}$ has the following form $$
\widehat{(\ell,\x)}=\sum_{k=0}^{N+1}\iprl{\lk,\hxk}=
\sum_{k=0}^N\iprl{\Qkw\Hk\pk,\yk}.
$$ The minimax a posteriori error has the following form $$
\hat{\sigma}=\bigl[1-\sum_{k=0}^N\iprl{\yk-\Hk\hxk,\Qkw\yk}\bigr]^\frac 12
\bigl(\sum_{k=0}^{N+1}\iprl{\lk,\pk}\bigr)^\frac 12,
$$ where $\pk$ solves
\begin{equation*}%   \label{eq:hzpkr}
   \begin{split}
     &\Fkt\hzk=\Ckt\hzko-\Hkt\Qkw\Hk\pk+\lk,
     \quad F'_{\ssstyle N+1}\hz_{\ssstyle N+1}=\ell_{\ssstyle N+1},\\
     &\Fko\pko=\Ck\pk+B_k\upQ_{\ssstyle 1,k}^{\ssstyle -1}B'_k\hzk,
     \quad{F}_0\p_0=S\iQz S'\hz_0,\,k=\overline{0,N}.
   \end{split}
\end{equation*}
\end{sz_thm}
Consider the minimax a posteriori estimates in the form of filters. 
\begin{sz_thm}\label{tm:fltr:r}
  Let $B_k=E$ and assume that the columns of the block matrix  $\left[\begin{smallmatrix}  
      F_k\\H_k\end{smallmatrix}\right]$
  are linear independent for any $k=\overline{0,N}$. Then for any $
  \ell\in\Rn$ the minimax a posteriori estimate of the linear function $\iprn{\ell,x_N}$ by observations $y_0\dots y_N$ in the form~\eqref{eq:yk:dae} has the following form 
  $
  \widehat{\iprn{\ell,x_N}}=\iprn{\ell,\hat{x}_{N|N}},
  $ where $\hx_{k|k}$ can be computed using the following algorithm: 
  \begin{equation*}%\label{eq:fltr:r}
    \begin{split}
      &\hx_{k|k}=P_{k|k}F'_k(Q^{-1}_{1,k-1}+C_{k-1}P_{k-1|k-1}C'_{k-1})^{-1}
      C_{k-1}\hx_{k-1|k-1}+P_{k|k}H'_k Q_{2,k}\yk,\\
      &P_{k|k}=\bigl(\Fkt(Q^{-1}_{1,k-1}+C_{k-1}P_{k-1|k-1}C'_{k-1})^{-1}\Fk+
      \Hkt\Qkw\Hk\bigr)^{-1},\\
      &P_{0|0}=(F'_0 \Qz F_0+H'_0 Q_{2,0}H_0)^{-1},\quad
      \hx_{0|0}=P_{0|0}H'_0 Q_{2,0}y_0,
    \end{split}
  \end{equation*}
\end{sz_thm}

Let us consider the contents of the third chapter in brief. 
Assume that the state $x(t)$ verifies the following DAE 
\begin{equation}
  \label{eq:Fx}
  \begin{split}
    &\ddt\Fx(t)-\Ct\xt=\ft,\\
    %&\Fx(a)=\fz,
  \end{split}
\end{equation}
and 
\begin{equation}
  \label{eq:Fxa}
  \Fx(a)=\fz
\end{equation}
where $F$ is a $m\times n$-matrix, $\Ct$ is a $
m\times n$-matrix with continuous on $[a,c]$ elements, $t\mapsto\ft\in\Rm$ is a vector-valued function from \Lm. In order to define the solution to~\eqref{eq:Fx} let us 
define a linear mapping $
\x\mapsto\oD\x\in\Lm\times\Rm$ by the following 
rule $$
  \begin{aligned}
    &\ov(\oD)\dfn\{\x\in\Ln:\Fx\in\Wm\}\dfn\WnF,\\
    &\oD\x\dfn(\ddt\Fx(t)-\Ct\xt,\Fx(a)),\x\in\ov(\oD).
  \end{aligned}
$$
%Denote by $\oDs$ the adjoint operator for $\oD$. 
Let $\tf\dfn(\f,\fz)\in\Lm\times\Rm$. 
Then $\xt$ is a solution of~\eqref{eq:Fx}-\eqref{eq:Fxa} if 
\begin{equation}
   \label{eq:Dx:d}
   \oD\x=\tf\,.
\end{equation}
It was proved in the dissertation that $\oD$ is 
closed dense defined linear mapping and its adjoit $\oDs$ was calculated. 
The minimax state estimation theory 
for linear DAEs in the form~\eqref{eq:Fx} can be constructed 
applying the same ideas as for the linear algebraic equations~\eqref{eq:Fx:a} 
presented above to operator equation~\eqref{eq:Dx:d}. This approach will be presented
below. 
   
Assume that $\x$ solves~\eqref{eq:Fx} and observed 
data \mbox{$t\mapsto\yt\in\Rl$} on $[a,c]$ is 
represented by
\begin{equation}  \label{eq:y}
  \yt=\Ht\xt+\et,
\end{equation}
where $\Ht$ is continuous $l\times n$-matrix on $[a,c]$, 
$t\mapsto\et\in\Rl$ is a realization of $l$-vector valued random process
with zero mean and continuous correlation function $
\R_\eta(t,s)=E\eta(t)\eta'(s)$. 
Define a linear mapping $\oH$ by the rule 
$\oH\xt = H(t)x(t)$.  

Let us consider a priori minimax estimates. 
We assume that the initial condition $
\fz\in\Rm$, input $
\f\in\Lm$ and correlation function $(t,s)\mapsto\R_\eta(t,s)$ 
are uncertain and belong to the given bounded set, that is: $$\tf\dfn[\f,\fz]\in\G,\R_\eta\in\Gw$$
As above we will look for the estimate of the linear transformation of the solution $\xt$ of~\eqref{eq:Fx}:  
\begin{equation*}
  \lx\dfn\int_a^c\iprn{\ell(t),\xt}\dt,\quad\ell\in\Ln
\end{equation*}%\label{fl:lx}
by means of the linear function $u$ of observed data
\begin{equation}
  \label{fl:uy}
  \uy\dfn\int_a^c\iprl{\ut,\yt}\dt+c,\quad \u\in\Ll,\,c\in\mathbb{R}.
\end{equation} 
We will refer the function of observation in the form~\eqref{fl:uy} as an 
estimate. Let us assign a worst-case estimation error  
\begin{equation*}%  \label{grnterr}
  \sigma(\u,c)\dfn\sup_{\x,\R_\eta}\{\M[\lx-\uy]^2|\oD\x\in\G,\R_\eta\in\Gw\}
\end{equation*}
The worst-case estimation error measures the quality of the estimate $\u$ 
and it does not depend on the particular realization 
of the uncertain parameters ($\fz,\f,\R_\eta$). 
\begin{sz_dfn}
The estimate $\hu_{\hat{c}}$ verifying 
  \begin{equation*}
    %\label{df:mmxerrnr:n}
     \sigma(\hu,\hat{c})\le
     \sup_{\x,\R_\eta}\{E[\lx-\uy]^2|\oD\x\in\G,\R_\eta\in\Gw\}=\sigma(u,c),
     \quad\u\in\setU_l,c\in\mathbb{R} 
  \end{equation*}
is called a minimax a priori mean-squared estimate. 
The number $
  \hat{\sigma}\dfn\inf_{\u,c}\sigma(\u,c)$ 
is called a minimax a priori mean-squared error. 
\end{sz_dfn}
The following propositions present the algorithms for calculation of a priori estimates
and errors in the variational form for the case of ellipsoidal bounding sets. 
\begin{sz_thm}\label{tm:mmxesterr:d}
  Let \begin{equation*}%  \label{setG}
    \G\dfn\{[\f,\fz]:\iprm{\Qz\fz,\fz}+\int_a^c\iprm{\Qn\f,\f}\dt\le 1\},
    \Gw\dfn\{\upR_\eta:\int_a^c\tr\bigl(\Qw\upR_\eta)\dt\le1\}
  \end{equation*}
  where $\upQ_{1,2}(t)$ are symmetric positive definite $m\times m$-matrices,  
 $\upQ^{-1}_{1,2}(t)$ are continuous on $[a,c]$, $\Qz$ is symmetric 
 positive definite $m\times m$-matrix. 
 Define a linear operator $\oT^+$ by the rule: $
\x\mapsto\oT^+\x=\hw$, where $\hw$ is a unique solution of the following 
optimization problem $$
\begin{aligned}
  &\|(\z_0,\z,\u)\|^2_{Q^{-1}}=
  \iprm{\iQz\z_0,\z_0}+\int_a^c\iprm{\iQn\z,\z}+\iprl{\iQw\u,\u}\dt\to
  \inf_{\z_0,\z,\u},\\
  &(\z_0,\z,\u)\in\setW\dfn
  \{\oDs(\z_0,\z)+\oHs\u=\x\}
\end{aligned}
$$ 
Then the minimax a priori mean-squared estimate is given by $$
  \huy=\ip{\oT^+\ell,(0,\y)}{L_2^m\times L_2^l}=\ipll{\hu,\y},
  $$ 
and minimax a priori mean-squared error may be represented as $$
  \hat{\sigma}=\sigma(\hu)=\|\oT^+\ell\|^2_{Q^{-1}}
$$ provided $$
  \ell\in\imT=\{\oDs(\z_0,\z)+\oHs\u,(\z_0,\z)\in\ov(\oDs),\u\in\Ll\}
  $$ 
\end{sz_thm}
\begin{sz_col}\label{cl:opstofunks:d}
  If the set $\imT$ is closed, then  
  \begin{equation}\label{eq:est:d}
  \hat{u}(t)=\Qw(t)\Ht\pt,\quad
  \sigma(\hu)=\int_a^c(\ell(t),\pt)_{\upR^n}\dt,
  \end{equation} where $\p$ solves the following two-point boundary value problem
  \begin{equation}\label{eq:ddtFpFz:d}\begin{aligned}
      &\ddt\Fz(t)=-\Ctt\zt+\Htt\Qw(t)\Ht\pt-\ell(t),\,\Fz(c)=0\\
      &\ddt F\pt=\Ct\pt+\iQn(t)\zt,\,
      F\p(a)=\iQz(FF^+\z(a)+d),F'd=0
  \end{aligned}
\end{equation}
\end{sz_col}
The next proposition represents a way to approximate the a priori estimate $\hat u$ by
means of Tikhonov regularization approach.
\begin{sz_thm}\label{tm:Ta:d}
  Take $\alpha_k>0$ and let $\p_k,\hz_k,d_k$ denote a unique solution of the following two-point boundary value problem: 
  \begin{equation}\label{GmltDscrEq}
    \begin{split}
      &\ddt\mA\pt=\Ct\pt+\alpha_k\iQn(t)\hzt,\,
      \upF\p(a)=\alpha_k\iQz(FF^+\hz(a)+d),\\
      &\ddt\Fz(t)=-\Ctt\zt+(\upE+\frac{1}{\alpha_k}\Htt\Qw(t)\Ht)\pt-
      \ell(t),\,\Fz(c)=0,
    \end{split}
  \end{equation}
  Then $$\begin{aligned}&\tilde{\ell}\in\imT
    \tilt
    \nrLl{\u_k-\hu}+\nrLm{\hz_k-\hz}+\|FF^+\hz_k(a)+d_k-\hz_0\|^2_{\upR^m}
  \xrightarrow{\alpha_k\downarrow0}0,
  \end{aligned}
  $$ where $\hz_0=FF^+\hz(a)+\hat{d}$, $\hu_k=\frac 1{\alpha_k}\Qw\upH\p_{k}$. 
% and $\hu$ denotes the a priori mean-squared minimax estimate. 
% of the linear function $
%   \tilde{\ell}=\mathrm{Pr}_{\cl\imT}(\ell)$. 
\end{sz_thm}
Let us present minimax estimates in the form of filters. 
\begin{sz_thm}\label{tm:fltr:d}
  Assume that $t\mapsto\Kt$ solves the following descriptor 
Riccati equation 
\begin{equation*}
  %\label{eq:rik1}
  \begin{split}
    &\ddt(\upF\Kt)=\Ct\Kt+\Ktt\Ctt-\Ktt\Htt\Qw\Ht\Kt+\iQn,\\
    &\upF\mathrm{K}(a)=FF^+\iQz FF^+,
  \end{split}
\end{equation*}
on $[a,c]$ and $t\mapsto\hzt$ verify the following differential-algebraic equation 
\begin{equation*}%\label{eq:obker}
  \ddt\Fz(t)+\Ctt\zt=\Htt\Qw(t)\Ht\Kt\zt,\quad\Fz(c)=\ell_0.
\end{equation*} 
The a priori minimax mean-squared estimate may be 
represented by 
$$
\widehat{(\ell,x(c))}=\int_a^c\iprl{\Qw(t)\Ht\Kt\hzt,\yt},
$$ and the a priori minimax mean-squared error is given by 
$$
\hat{\sigma}=\iprm{FK(c){F^+}'\ell_0,{F^+}'\ell_0}.
$$
Let $\hx$ denote a solution of the following 
linear differential-algebraic equation 
\begin{equation*}%  \label{eq:fltr}
  \begin{split}
    &\ddt\upF\hxt=\Ct\hxt+\Ktt\Htt\Qw(t)(\yt-\Ht\hxt),\\
    &\upF\hx(a)=0
  \end{split}
\end{equation*}
Then the a priori minimax mean-squared estimate is given by $$
    \widehat{(\ell_0,x(c))}=\iprm{\upF\hx(c),{F^+}'\ell_0}.
$$
\end{sz_thm}
Definitions and representations for the minimax a posteriori estimates are given in the dissertation. 

{ \bf Conclusion. } The dissertation presents a generalized Kalman duality concept and minimax 
state estimation approach for linear differential-algebraic equations. The key results of the dissertation are as follows: 
\begin{itemize}
\item generalized Kalman duality concept;
\item new definitions of the minimax estimates based on the generalized Kalman duality concept;
\item new variational form of minimax state estimation algorithms 
for linear DAE;
\item new minimax state estimation algorithms for DAE in the form of filters;
\item efficient description of the reachability set for DAE;
\item description of the uncertainty propagation by the DAE;
\end{itemize}

{  \bf List of papers and conference presentations related to the dissertation. }
%\bibliographystyle{unsrt}
%\bibliographystyle{alpha}
%bib sources for bib-tex
\nobibliography{refs/myrefs,refs/refs}
\paragraph{International journals (peer-reviewed)}
\mybeglist
\mybibentry{Zhuk2004}
\myendlist
\paragraph{International journals}
\mybeglist
\mybibentry{Zhuk2005a}
\myendlist
\paragraph{National journals (peer-reviewed)}
\mybeglist
\mybibentry{Zhuk2005b}
\mybibentry{Zhuk2005}
\myendlist
\paragraph{National journals}
\mybeglist
\mybibentry{Zhuk2004a}
\mybibentry{DolenkoZhuk2002}
\myendlist
\paragraph{Preprints}
\mybeglist
\mybibentry{Zhuk2006c}
\mybibentry{Zhuk2006b}
\myendlist
%\newpage
\noindent\paragraph{International conferences}
\mybeglist
\mybibentry{Zhuk2006a}
\mybibentry{Zhuk2006}
\mybibentry{Zhuk2005c}
\mybibentry{Zhuk2004b}
\mybibentry{Zhuk2003}
\myendlist
\paragraph{National conferences}
\mybeglist
\mybibentry{Zhuk2004c}
\myendlist
\paragraph{References}
\mybeglist
\mybibentry{Tikhonov1977}
\mybibentry{Hanke1989}
\mybibentry{Dai1989}
\mybibentry{Bertsekas1971}
\mybibentry{Tempo1985}
\mybibentry{Kuntsevich1992}
\mybibentry{Chernousko1994}
\mybibentry{Kurzhanski1997}
\mybibentry{Nakonechny2004}
\myendlist

\end{document}